\documentclass{amsart}

\usepackage{amssymb, amsmath, latexsym, a4}
\usepackage[latin1]{inputenc}
\usepackage[all]{xy}

\usepackage{tikz}
\usepackage{pinlabel}
\usepackage{amsfonts}
\usepackage{amscd}

\usetikzlibrary{matrix,calc,arrows}
\newcommand*\Z{\mathbb{Z}}

\newcommand\xto[1]{\xrightarrow[]{#1}}

\newcommand{\pn} {\par \noindent}

\newcommand{\pbn}{\par \bigskip \noindent}

\newtheorem{De}{Definition}[section]

\newtheorem{Prop}[De]{Proposition}
\newtheorem{Le}[De]{Lemma}
\newtheorem{Co}[De]{Corollary}

\newcommand{\rdg}{\hfill $\Box $}

\def\t{\otimes }

\def \HL{{\sf HL}}
\def \g{{\bf g}}
\def \h{{\bf h}}

\def \Li{{\sf Lie}}

\begin{document}
\def\t{\otimes }
\def\al{\alpha}
\def\bb{\beta}
\def\dd{\delta}
\def\te{\theta}

\title{A subcomplex of Leibniz complex}

\author[Teimuraz Pirashvili]{Teimuraz  Pirashvili}
\address{
Department of Mathematics\\
University of Leicester\\
University Road\\
Leicester\\
LE1 7RH, UK}

\maketitle

\begin{abstract} Using the free graded Lie algebras we introduce a natural subcomlex of the Loday's complex of a Leibniz algebra. \end{abstract}

\bigskip
\section{Introduction}

This note, except for the introduction and the second half of Section \ref{conj}, were written in 2002 during my stay in Nantes, where I was invited by Vincent Franjou.

Let $\h$ be a Leibniz algebra and let ${\sf CL}_*(\h)$ be the chain compex constructed by Loday \cite{HC}. The aim is to intriduce a natural subcomplex $L(\h,1)_*$ of  ${\sf CL}_*(\h)$, which conjecturaly computes the Quillen derived functors of the
left adjoint functor to the inclusion
$$ {\sf Lie\ Algebras} \subset {\sf Leibniz \ Algebras}.$$

The paper also includes a conjecture of Jean-Luis Loday (see Section \ref{conj}).

We refer the reader to \cite{perfect} about other conjectures on Leibniz homology.

\section{Definition of the subcomplex}

\pbn We are working over a field $K$ of characteristic zero.

\pn A Leibniz algebra $\h$ is a vector space equipped with a bilinear operation
$$\{-,-\}:\h\t\h\to \h$$
satisfying the Leibniz identity
$$\{x,\{y,z\}\}=\{\{x,y\},z\}-\{\{x,z\},y\}.$$
It is clear that Lie algebras are Leibniz algebras for which
the antisymmetry condition $\{x,y\}+\{y,x\}=0$ holds.

We refer the reader to \cite{LP}, \cite{HL} and \cite{Lrep}  for more on Leibniz algebras, Leibniz homology and  Leibniz representations.

 For any Leibniz algebra
$\h$ we let $\h_{\sf Lie}$ be the quotient of $\h$ by the relation 
$\{x,y\}+\{y,x\}=0$. In this way one obtains the functor
$$(-)_{\sf Lie}:{\sf Leibniz\  algebras}\to {\sf Lie \ algebras}$$
which is the left adjoint of the inclusion functor 
$ {\sf Lie \ algebras}\subset {\sf Leibniz\  algebras}.$ 

For Leibniz algebras we have {\it Leibniz homology} introduced by Jean-Louis Loday 
in \cite{HC}, see also \cite{LP}. For a Lie algebra $\h$, the Leibniz
 homology (denoted by  $\HL_*(\h)$) is defined as the 
homology of the complex ${\sf CL}_*(\h)$. By 
definition ${\sf CL}_n(\h)=\h^{\t n}$, while the Loday boundary 
is given by
$$d(x_1\t\cdots \t x_n)=\sum_{i<j}(-1)^{j+1}(x_1\t \cdots \t x_{i-1}\t[x_i,x_j]\t\cdots \t \hat{x_j}\t\cdots \t x_n).$$
Thus 
$${\sf CL}_*(\h)= \ \ (\cdots \to \h^{\t n}\to \h^{\t n-1}\to \cdots \to \h^{\t 2}\to \h \to K).$$
The boundary map $\h\to K$ in the complex is the zero map.

\
\pn Let us recall that a graded 
Lie algebra is a graded vector space 
$\g:= \bigoplus _{n\in \Z}\g _n$ equipped with a bilinear map
$$[-,-]:\g _n\t \g _m\to \g _{n+m}, \ n,m\in \Z$$
satisfying the following equations
$$[x,y]+(-1)^{\mid x\mid \mid y\mid}[y,x]=0,$$
$$[x,[y,z]]=[[x,y],z]-(-1)^{\mid y\mid \mid z\mid}[[x,z],y].$$
Here $x,y,z$ denotes homogeneous elements and we assume that
 $x\in \g _{\mid x\mid}$.
If $A_*$ is a graded associative algebra then it can be considered also as 
a graded Lie algebra by $[x,y]:=xy-(-1)^{\mid x\mid \mid y\mid}yx$.

\ 

Let $V$ be a vector space. We let $T(V,1)$ denote the tensor
algebra on $V$, where $V$ is concentrated in the degree 1. In other words
$$T(V,1)=\bigoplus _{n\geq 0} T(V,1)_n, \ \ T(V,1)_n=V^{\t n}, n\geq 0.$$
We let $L(V,1)=\bigoplus _{n\geq 0} L(V,1)_n, $ be the graded 
Lie subalgebra of $T(V,1)$ generated
by $V$.  The Lie algebra $L(V,1)$ is called {\it free graded Lie algebra
generated by $V$ concentrated in the degree 1}. It is clear that
$L(V,1)_1\cong V$ and the map $x\t y\mapsto [x,y]$ yields the isomorphism
$S^2V\cong L(V,1)_2,$ where $S^2$ denotes the second symmetric power. 
The vector space
$L(V,1)_n$ is spanned by the elements of the form 
$[x_1,[x_2,[\cdots [x_{n-1} ,x_n]\cdots]$ with $x_1,\cdots, x_n\in V$. 
For simplicity we denote this type of element by $[[x_1,\cdots, x_n]]$. 

\

We let $i_n:L(V,1)_n\subset T(V,1)_n=V^{\t n}$ be the inclusion. The main
observation of this note  is the fact that if $\h$ is a Leibniz algebra, 
then $L(\h,1)_*$ is closed under Loday boundary map, and therefore $L(\h,1)_*$ is a subcomplex of ${\sf CL}_*(\h)$. In other words 
one has the following
 
\begin{Prop}\label{qvekomplexi} Let $\h$ be a Leibniz algebra. 
Then  $$d(i_n(L(\h,1)_n))\subset i_{n-1}(L(\h,1)_{n-1}).$$
\end{Prop}

In order to describe the induced boundary map 
$\dd:L(\h,1)_n\to L(\h,1)_{n-1}$, we need the linear map 
$p_n:T(V,1)_n\to L(V,1)_n$ given by
$$p_n(x_1\t\cdots  \t x_n)=[[x_1,\cdots ,x_n]].$$

\begin{Prop}\label{komplexi} 
Let $\h$ be a Leibniz algebra. Then, for any element 
$\omega\in L(\h,1)_{n+1}$ one has
 $$di_{n+1}(\omega)=(-1)^np_n\circ f^n\circ i_{n+1}(\omega),$$ 
where
$$f^n: \h^{\t n+1}\to \h^{\t n}$$
is the linear map given by
$$f^n(x_1\t\cdots \t x_n\t x_{n+1})=
x_1\t \cdots \t x_{n-1}\t \{x_n,x_{n+1}\}.$$

\end{Prop}

\bigskip

We give the proof of Proposition \ref{qvekomplexi} and   Proposition \ref{komplexi} in Section \ref{proof1}.

\pbn As a corollary we get for any Leibniz algebra $\g$ a chain complex 
$$L(\h,1)_*:= (
\cdots \to L(\h,1)_{n+1}\xto{\delta_n} L(\h,1)_{n}\xto{\delta_{n-1}}  \cdots \
\xto{\delta_2} L(\h,1)_{2}\xto{\delta_1}
L(\h,1)_{1})$$
The first nontrivial boundary maps are given by
\begin{align*}
\delta_1([x,y])&=\{x,y\}+\{y,x\},\\
\delta_2([x,[y,z]])&=[x,\{y,z\}]+[x,\{z,y\}]-[y,\{z,x\}]-[z,\{y,x\}]\\
\delta_3([x,[y,[z,u]]])&=[x,[y,\{z,u\}]]+[x,[y,\{u,z \}]]-[x,[u,\{z,y\}]] -[x,[z,\{u,y\}]]\\
&+ [y,[u,\{z,x\}]]+[y,[z,\{u,x\}]]-[z,[u,\{y,x\}]] -[u,[z,\{y,x\}]]
\end{align*}
We let $\Li_n(\h)$, $n\geq 1$ be the  homology of the complex 
$L(\h,1)_*$.

\begin{Le} For any Leibniz algebra $\h$ the group $\Li_1(\h)$ is 
isomorphic to the Liezation of $\h:$
$$\Li_1(\h)\cong \h_{\sf Lie}$$ 
\end{Le}

\section{Two conjectures}\label{conj}

The fisrt one says:
 
{\bf Conjecture 1}. If  $\h$ is a free Leibniz algebra, then $\Li_n(\h)=0$ for $n>1$.

If this conjecture is true, then  $\Li_*(\h)$ computes 
the Quillen derived functors of the functor
$$(-)_{\sf Lie}:{\sf Leibniz \ Algebras}\to {\sf Lie\ Algebras}.$$

Since 2002, I have had the oportunity to discuss the conjecture with several people, including Jean-Luis Loday, Michael Robinson and Christine Vespa, but unfortunately, our attemts to prove the conjecture failed. 

The second one claims that our chain complex is only a part of a big pichture:

 {\bf Conjecture 2}. (J.-L. Loday. 2006 \cite{0707} ). Let $e_n^{(i)}$ be the image of the $i$-th Eulerian idempotent on $\g^ {\otimes n}$. Then
$$d(e_n^{(i)}) \subset  e_{n-1}^{(1)}\oplus \cdots \oplus e_{n-1}^{(i)}.$$
 For definition of the Eulerian idempotents see \cite[pp 142-144]{HC}.

%Since the idempotent $e^1$ gives the free Lie algebra, the main result of this notes shows that the conjecture 2 is true %for $i=1$. In the same email he had  wrote  that for $i=2$, he checked his conjecture for small values of $n$. 
%As a consequence of this conjecture there is a filtration on the  
%Leibniz complex given by
%$$F_i =  e_{n}^{(1)}\oplus \cdots \oplus e_{n}^{(i)}$$

\section{A generalization of  Wigner's identity}
In this section $V$ denotes a graded vector space and $\t$ denotes the tensor
product of graded vector spaces.
Any  linear map $D:V\to V$ has an extension  as a linear map $D_n:T(V)\to T(V)$
given by 
$$D_n(v_1\t \cdots 
\t v_n)=\sum_
{i=1}^n v_1\t \cdots \t D(v_i)\t \cdots \t v_n.$$
In particular $D_0(1)=0$.  Moreover $D_n(\omega)=n\omega$, provided $D=Id$.
Here $\omega\in T^n(V).$

Let us recall that $T(V)$ is a graded Hopf algebra with following coproduct:
$$\Delta (v)=v\t 1+1\t v, \,v\in V.$$
We let $\mu$ be the multiplication map $T(V)\t T(V)\to T(V)$ and 
we let $L(V)$ be the Lie  subalgebra of $T(V)$ generated by $V$. Let
$$p_D:T^n(V)\to L(V)_n$$
be the map given by 
$$ p_D(v_1\t \cdots \t v_n)=[v_1,[v_2,\cdots ,[v_{n -1}, Dv_n]\cdots ].$$

\begin{Prop}\label{wigner} For any 
$\omega\in T^n(V)$ one has $$\mu\circ (p_D\t Id_{T(V)})\circ \Delta (\omega)=D_n(\omega )$$
\end{Prop}

Let us observe that Proposition \ref{wigner} for  $D=Id$ is  a  result 
of \cite{wigner}. The proof in our setup is completely similar to one
given in {\it loc.cit.}

{\it Proof}. The result is obvious for $n=1$ and therefore we 
can use the induction on $n$. For an  element $\omega$ of $T^n(V)$ we use the notation $\omega = v_1 \cdots v_n$, $v_1,\cdots, v_n\in V$. We may write
$\Delta (v_2\cdots v_n)=1\t v_2\cdots v_n+\sum _i a_i\t b_i$ 
with each 
$a_i\in  \bigoplus _{i\geq 1}T^i(V)$. The induction hypothesis gives 
$\sum_{i} p_D(a_i)b_i=\sum_{k\geq 2}v_2\cdots D(v_k)\cdots v_n.$ On the other hand, we have $p_D(v\tau)=[v,p_D(\tau)]=vp_D(\tau)-
(-1)^{\mid \tau \mid}p_D(\tau)v.$ We have also
$$\Delta(\omega)=(v_1\t 1+1\t v_1)(1\t v_2\cdots v_n+\sum _ia_i\t b_i)$$
$$= v_1\t v_2\cdots v_n+ \sum_iv_1a_i\t b_i+1\t \omega +(-1)^{\mid a_i\mid}
\sum_ia_i\t v_1b_i.$$
Thus $$\mu(p_D\t Id)\Delta (\omega)= (Dv_1)v_2\cdots v_n+\sum_ip_D(v_1a_i)b_i+
+(-1)^{\mid a_i\mid}\sum_i p_D(a_i)v_1b_i$$
$$= (Dv_1)v_2\cdots v_n+\sum_iv_1p_D(a_i)b_i$$
$$=(Dv_1)v_2\cdots v_n+ v_1\sum_{k\geq 2}v_2\cdots D(v_k)\cdots v_n$$
\rdg

\begin{Co}\label{friedrichs} For any $\omega\in L(V)_n$ one has
$$p_D(\omega)=D_n(\omega ).$$
\end{Co}

{\it Proof}. If $\omega \in L(V)_n$ then 
$\Delta(\omega)=\omega \t 1+1\t \omega$. 
Thus the result is a direct consequence of Proposition \ref{wigner}.\rdg

\section{Proof of Propositions 
\ref{qvekomplexi} and \ref{komplexi}}\label{proof1}
In this section we let 
$\dd=p_n\circ f^n\circ i_{n+1}:L(\h,1)_{n+1}\to L(\h,1)_n$. 
We will prove the following

\begin{Prop}\label{orive}
$di_{n+1}= (-1)^ni_n\dd.$
\end{Prop}
Since $i_k$ is injective for all $k$, this implies both Propositions 
\ref{qvekomplexi} and \ref{komplexi}.

We start with general remarks.
Let $\g$ be a Lie algebra and let $V$ be a left $\g$-module. For $x\in V$
and $g\in \g$ we let $\{x,g\}$ be the result of the action of $g$ on $x$. It
is well known that $\g$ acts also on $V^{\t n}$ by 
$$\{x_1\t\cdots \t x_n,g\}:=\sum_{i=1}^n x_1\t\cdots \t \{x_i,g\}\t\cdots \t x_n.$$

\begin{Le} One has the following identity:
$$\{[[x_1,\cdots,x_n]], g\}=\sum_{i=1}^n[[x_1,\cdots,\{x_i,g\},\cdots ,x_n]]$$
\end{Le}

{\it Proof}. For $n=2$  we have
$$\{[x_1,x_2],g\}=\{x_1\t x_2+x_2\t x_1,g\}$$
$$=\{x_1,g\}\t x_2+x_1\t \{x_2,g\}+\{x_2,g\}\t x_1+x_2\t \{x_1,g\}$$
$$=[\{x_1,g\} ,x_2]+[x_1,\{x_2,g\}].$$
For $n\geq 3$ we have $[[x_1,\cdots,x_n]]=[x_1, [[x_2,\cdots,x_n]]$ and Lemma
follows by the induction.

\begin{Co} For any $n\geq 1$ the vector space $L(V,1)_n$ is a $\g$-submodule 
of $T(V,1)_n=V^{\t n}$ and the map $p_n$ is $\g$-linear.  
\end{Co}

Let $\h$ be a Leibniz algebra. Since $\{ x,\{y,z\}+\{y,z\}\}=0$, it follows
that the Lie algebra $\g:= \h_{\sf Lie}$ acts on $h$ via $\{x,{\bar y}\}:=
\{x,y\}$. Here $\bar y$ denotes the class of $y\in \h$ in $\g$. Let us observe
that $f(x,y)=\{x,y\}$ defies a $\g$-linear map $f:\h\t \h\to \h$. Thus $f^n$ is
also a $\g$-linear map and as a consequence $\dd=p_n\circ f^n\circ i_{n+1}$ is aslo a
$\g$-linear. 

\begin{Le}\label{fr} For any $x\in \h$ and $\omega \in L(\h,1)_{n-1}$ one has
$$\delta ([x,\omega ] )=[x,\delta \omega]- (-1)^{\mid \omega \mid}\{\omega, {\bar x}\}$$
\end{Le}

{\it Proof}. We have 
$$\delta ([x,\omega ] )-[x,\delta \omega]=-(-1)^{\mid \omega \mid}p_{n-1}f^{n-1}(i_{n-1}\omega \t x).$$
We let $D:\h\to \h$ be the map given by $D(v)=\{v,x\}$. Since $f^{n-1}=Id_{{\h}^{\t n-2}}\t f$, we have
$$p_{n-1}(Id_{h^{\t n-2}}\t f)(i_{n-1}\omega \t x)=p_D(\omega)$$ and
Lemma is a direct consequence of Corollary \ref{friedrichs}. \rdg

For a linear map $D:V\to V$ we let $\tilde{D}_n:V^{\t n}\to V^{\t n}$ be the 
map
given by 
$${\tilde D}_n(v_1\t \cdots 
\t v_n)=\sum_
{i=1}^n (-1)^i D(v_i)\t v_1\t \cdots \t \hat{v}_i\t \cdots \t v_n.$$

\begin{Le}\label{lis} 
For any $\omega\in L(V,1)_n$, $n\geq 2$ and any linear map 
$D:V\to V$ one has
$${\tilde D}_n(\omega)=0$$
\end{Le}

{\it Proof}. We work by induction. If $n=2$ we have
$${\tilde D}_2(x\t y+y\t x)=-D(x)\t y+D(y)\t x -D(y)\t x+D(x)\t y=0.$$
If $n>2$ then we may assume that $\omega =[x,\tau]$, where $\tau \in 
L(V,1)_{n-1}.$ Thus $\omega=x\t \tau - (-1)^{\mid \tau \mid}\tau \t x$ and 
therefore we have
$$\tilde {D}_n(\omega)=-D(x)\t \tau- j_2(x)(\tilde {D}_{n-1}(\tau)) - 
(-1)^{\mid \tau \mid }\tilde {D}_{n-1}(\tau) \t x +D(x)\t \tau=0$$
thanks to the induction assumption. Here for any $x\in V$ we let $j_2(x):V^{\t n-1}\to V^{\t n}$ be the map given by
$y_1\t\cdots y_{n-1}\mapsto \t y_1\t x\t \cdots \t y_{n-1}$. \rdg

The following lemma is an immediate consequence of the definitions (compare 
also with formula  (10.6.3.1) of \cite{HC}).

\begin{Le} For any $x\in \h$ and $\omega \in V^{\t n}$ one has
$$d (\omega \t x )= d(\omega)\t x +(-1)^{\mid \omega \mid}\{\omega, \bar{x}\}$$
and
$$d (x\t \omega )= -x\t d(\omega) + \tilde{D}_{\mid \omega \mid}(\omega), $$
for $D=\{x,-\}:\h\to \h$.
\end{Le}

\begin{Le}\label{fro} For any $x\in \h$ and $\omega \in L(\h,1)_n$ one has
$$d ([x,\omega ] )=-[x,d\omega]- \{\omega,{\bar x}\}.$$
\end{Le}

{\it Proof}. Since $[x,\omega]= x\t \omega -(-1)^n\omega \t x$, we have
$$d ([x,\omega ] )=d(x\t \omega -(-1)^n\omega \t x)$$
By the previous Lemma, this expresion can be rewriten:
$$-x\t d\omega +\tilde{D}_n(\omega)-(-1)^nd\omega\t x-(-1)^n(-1)^n\{\omega, \bar{x}\}$$
Thanks to Lemma \ref{lis} $\tilde{D}(\omega)=0$ and we get the result.\rdg

 Now we are in position to prove  Proposition \ref{orive}.

{\it Proof  Proposition \ref{orive}}. We use the induction on $n$. First consider the case $n=1$. So we 
may assume that $\omega= [x,y]=x\t y+ y\t x$. Hence
$$di_2(\omega)=d(x\t y+ y\t x)=-\{x,y\}-\{y,x\}=-\dd([x,y]).$$
If $n>1$, then we may assume that $\omega=[x,\tau]$, with $\tau \in L(V,1)_n$.
Thanks to Lemma \ref{fro} we have 
$$di_{n+1}(\omega)= d([x,\tau])= -[x,d\tau]- \{\tau,{\bar x}\}$$
and the induction assumption gives
$$di_{n+1}(\omega)= (-1)^n[x,\dd \tau]-\{\tau,{\bar x}\}$$
and the result follows from Lemma \ref{fr}.

\end{document}